\pdfoutput=1
\documentclass[12pt]{article}
\usepackage{amsmath}
\usepackage{amsfonts}
\usepackage{amssymb}
\usepackage{latexsym}
\usepackage[dvips]{graphicx}
\usepackage{float}
\usepackage{geometry}
\AtEndDocument{\bigskip{\footnotesize
\textsc{Department of Mathematics, University of Macau, Macao}\par
\textit{E-mail address: }\texttt{maiweixiong@gmail.com}\par
\medskip
\textsc{Department of Mathematics, University of Macau, Macao}\par
\textit{E-mail address: }\texttt{fsttq@umac.mo}

}}

\usepackage{lipsum}

\newtheorem{thm}{\bf Theorem}[section]
\newtheorem{lem}[thm]{\bf Lemma}
\newtheorem{prop}[thm]{\bf Proposition}

\newenvironment{proof}{\noindent{\em Proof:}}{\quad \hfill$\Box$\vspace{2ex}}

\newenvironment{remark}{\noindent{\bf Remark}}{\vspace{2ex}}

\begin{document}
 \title{Aveiro Method in Reproducing Kernel Hilbert Spaces Under Complete Dictionary\let\thefootnote\relax\footnotetext{This work was supported by University of Macau research grant MYRG116(Y1-L3)-FST13-QT and Macau Government FDCT 098/2012/A3.}}
\author{
Weixiong Mai,
and Tao Qian
}
\date{}
 \maketitle
\begin{abstract}
Aveiro Method is a sparse representation method in reproducing kernel Hilbert spaces (RKHS) that gives orthogonal projections in linear combinations of reproducing kernels over uniqueness sets. It, however, suffers from determination of uniqueness sets in the underlying RKHS. In fact, in general spaces, uniqueness sets are not easy to be identified, let alone the convergence speed aspect with Aveiro Method. To avoid those difficulties we propose an anew Aveiro Method based on a dictionary and the matching pursuit idea. What we do, in fact, are more: The new Aveiro method will be in relation to the recently proposed, the so called Pre-Orthogonal Greedy Algorithm (P-OGA) involving completion of a given dictionary. The new method is called Aveiro Method Under Complete Dictionary (AMUCD). The complete dictionary consists of all directional derivatives of the underlying reproducing kernels. We show that, under the boundary vanishing condition,  bring available for the classical  Hardy and Paley-Wiener spaces, the complete dictionary enables an efficient expansion of any given element in the Hilbert space. The proposed method reveals new  and advanced  aspects in both the Aveiro Method and the greedy algorithm.

\end{abstract}


\section{Introduction}
We first give a revision on the theory of reproducing kernel Hilbert spaces (e.g. \cite{NA,cfrs12, cfrst}).\\
Let $\mathcal H$ be a Hilbert space, and $E$ an abstract set and ${\bf h}$ a $\mathcal H$-valued function on $E$.
Then, consider the linear transformation
\begin{equation}\label{2.1}
f(p) = \langle {\bf f},{\bf h}(p)\rangle_{\mathcal H}, \quad {\bf f} \in \mathcal H\,,
\end{equation}
from $\mathcal H$ into the linear space ${\mathcal F}(E)$ comprising all the complex valued
functions on $E$. Contruct a positive definite quadratic form function
\begin{equation}\label{2.2}
 K(p,\overline q) = \langle {\bf h}(q) , {\bf h}(p)\rangle_{\mathcal H}~ \quad {\rm on}~\quad E {\rm \times} E.
\end{equation}

The following is the basic theory of reproducing kernel Hilbert spaces.

\begin{prop}\label{pro:1.1}
\begin{itemize}
\item[{\rm (I)}] \quad The range of the linear mapping $(\ref{2.1})$
by $\mathcal H$ is characterized
as the reproducing kernel Hilbert space $H_K(E)$ admitting the reproducing  kernel $K(p,\overline q)$ whose characterization is given by the two properties: $(i)$ $K(\cdot,\overline q)  \in H_K(E)$ for any $q \in E$ and, $(ii)$ for any $f \in H_K(E)$ and for any $p \in E$, $\langle f(\cdot),K(\cdot,\overline p)\rangle = f(p)$.
\item[{\rm (II)}]  \quad
In general, we have the inequality
$$\Vert f \Vert \leq \Vert {\bf f} \Vert_{\mathcal H}.$$
\noindent
Here, for any member $f$ of $H_K(E)$ there exists a uniquely determined $
{\bf f}^* \in \mathcal H$ satisfying
$$f(p) = \langle {\bf f}^*,{\bf h}(p)\rangle_{\mathcal H}~ \quad {\rm on} \quad E$$
and
\begin{equation}\label{2.3}
\Vert f \Vert = \Vert {\bf f}^* \Vert_{\mathcal H}.
\end{equation}
\item[{\rm (III)}] \quad
In general, we have the inversion formula in $(\ref{2.1})$ in the form
\begin{equation}\label{2.4}
 f \mapsto {\bf f}^*
\end{equation}
in $(II)$ by using the RKHS $H_K(E)$.
\end{itemize}
\end{prop}

\noindent For more information on reproducing kernel Hilbert space (RKHS), please see \cite{NA,Sa97,Sa10}.

In \cite{cfrs12, cfrst} S. Saitoh et al propose the so-called Aveiro Method aiming to  construct an approximating function of $\bf {f^{*}}$ involving a finite number of sampling points of $E$.
\begin{prop}\label{prop2.4}
Suppose that $\{p_j\}_{j=1}^n$ are $n$ distinct points in $E$. Define a Hermitian matrix $A_n$ with the elements
\begin{equation}\label{ajk}
a_{j,k} =  \langle {\bf h}(p_{k}),{\bf h}(p_j)\rangle_{\mathcal H},
\end{equation}
and further assume that $A_n$ is positive definite.

For \begin{equation}
f(p) = \langle {\bf f},{\bf h}(p)\rangle_{\mathcal H}, \quad {\bf f} \in \mathcal H\, ,
\end{equation}
we have
\begin{equation}\label{eq8}
{\bf f}^*_{A_n} = \sum_{j=1}^n \sum_{k=1}^n f(p_j) \widetilde{ a_{j,k}^{(n)}} {{\bf h}(p_{k})}
\end{equation}
satisfying
\begin{equation}\label{interpolation}
f(p_j) = \langle {\bf f}^*_{A_n},{\bf h}(p_j)\rangle_{\mathcal H},
\end{equation}
where $\widetilde{a^{(n)}_{j,k}}$'s are the elements of $\overline{A_n^{-1}}$
(here we use the notation $\widetilde {a_{j,k}^{(n)}}$ because elements of $\overline {A_n^{-1}}$ depend on $n$).
Moreover, if there exists any ${\bf g}\in\mathcal H$ such that
$$
f(p_j)=\langle {\bf g}, {\bf h}(p_j)\rangle_{\mathcal H},
$$
we have
\begin{equation*}
||{\bf f}_{A_n}^*||_{\mathcal H}\leq ||{\bf g}||_{\mathcal H}.
\end{equation*}
\end{prop}

\noindent The convergence of ${\bf f}^*_{A_n}$ is based on the following Proposition.
\begin{prop}\label{convergence}
 Let $\{p_j\}_{j=1}^\infty$ be a sequence of distinct points in $E$, that is of the positive definiteness property set as in Proposition 1.2 for any $n,$ and a uniqueness set for the RKHS $H_K$; that is, for any $f \in H_K$, if $f(p_j) = 0$ for all $p_j,$ then $f  \equiv 0$.
Then, in the space $\mathcal H$
\begin{equation}
\lim_{n \rightarrow \infty} {\bf f}^*_{A_n} =  {\bf f}^*
\end{equation}
for ${\bf f}^*$ given in Proposition 1.1.

\end{prop}

As shown in Proposition \ref{convergence}, the convergence of ${\bf f}^*_{A_n}$ depends on the assumption that $\{p_j\}^{\infty}_{j=1}$ is a uniqueness set. The effectiveness of Aveiro Method is not guaranteed in practical use due to the following reasons: (1) The knowledge of uniqueness sets of a RKHS is usually not sufficient, and uniqueness sets are not easy to be identified; and (2) For a uniqueness set $\{p_j\}^{\infty}_{j=1}$, there do not exist results addressing convergence behavior in terms of ${A_n}.$ What is known would be only that the series converges to ${\bf f}^\ast.$
 In this paper we propose an modified Aveiro Method over  a given dictionary, as well as the matching pursuit idea (e.g. \cite{MZ,Tem}). Owing to these features it becomes practical. By doing this we not only avoid the uniqueness set issue but also increase the convergence speed. More specifically, the proposed theory and algorithm depend on a recently proposed concept, the completion of the dictionary originally given (\cite{Q}). We call the proposed method Aveiro Method Under Complete Dictionary (AMUCD).

 Within the above axiomatic formulation of RKHS we will be working with the simple cases in which $\mathcal H=H_K(E),$ and thus $\mathcal H$ itself is a RKHS.
  In the original Aveiro Method, by the definition of the Aveiro representation, $f^*_{A_n}$ is meaningless if $p_j=p_k$ for some $j\neq k$. The proposed AMUCD, in a sense, allows the cases $p_j=p_k$ for $j\neq k.$ It is done by involving the complete dictionary consisting of the original dictionary elements together with all the possible directional derivatives, reducing to derivatives in particular cases, of the dictionary elements. In our case a dictionary consists of reproducing kernels. It is by introducing the complete dictionary concept that enables Qian to propose a new type of greedy algorithm called Pre-Orthogonal Greedy Algorithm (P-OGA) in \cite{Q}. It is shown that P-OGA is among the most effective matching pursuit methods. Applying the same idea, in AMUCD the matrices $A_n$ can involve derivatives of the reproducing kernels. In the meantime, as in AFD (see \cite{QWa1}, \cite{Q}), the representation and its derivatives are of interpolation property at the selected points. AMUCD is, in fact, an alternative representation of AFD. In such way the capacity of Aveiro method is considerably lifted up. Denote by $\widetilde f_{A_n}^*$ the revised $f_{A_n}^*$.
For $n$ points $\{p_j\}_{j=1}^n$ in $E$, we select the next point that satisfies
\begin{align}\label{min_proH}
p_{n+1}^*= \arg \min_{p_{n+1}\in E } ||{f}-\widetilde {f}^{*}_{A_{n+1}}||.
\end{align}
 By using a complete dictionary, the existence of $p_{n+1}^*$ is evident when the remainder energy approaches to zero at the boundary (see Lemma \ref{lem_new}). We call such property $\lq\lq$boundary vanishing condition (BVC)". If BVC holds and each element of $\{p_j\}_{j=1}^\infty$ is selected according to the principle (\ref{min_proH}), then $\{p_j\}_{j=1}^\infty$ does not have to be a uniqueness set of $\mathcal H$.  Thus the convergence of $\widetilde{f}_{A_n}^*$ is not a conclusion of Proposition 1.3, and, instead, requires a separate proof. We give a proof of the convergence
of such formulated $\widetilde {f}_{A_n}^*$ in \S 2.

 In \S 3 we perform AMUCD to two special cases, the Hardy space $H^2(\mathbb D)$ and the Paley-Wiener space $W(\frac{\pi}{h}), h>0$.  For $H^2(\mathbb D)$, we show that BVC holds (see also \cite{WQ}). For $W(\frac{\pi}{h})$ we are able to show a weak BVC property.

In this study we show that AMUCD is identical with P-OGA  (see \S 2). AMUCD has the advantage of not working out the related orthonormal system. As a matter of fact, in many cases P-OGA do not have explicit formulas for the related orthonormal system functions like the Takenaka-Malmquist (TM) system case in the classical Hardy spaces case.

The writing plan of the paper is as follows. In \S 2 we present AMUCD, giving a set of sufficient conditions for existence of $p_{n+1}^*$ and proving the convergence of $\widetilde f_{A_n}^*$ without assuming $\{p_j\}^\infty_{j=1}$ being a uniqueness set. In \S 3 we implement our method to the Hardy and Paley-Wiener spaces.

\medskip

\section{Aveiro Method Under Complete Dictionary}
In this section, we present a revised ${\bf f}^*_{A_n}$. Then, we present AMUCD. In particular, we are concerned with the cases $\mathcal H=H_K(E)$, where $E$ is a domain in $\mathbb C$ and ${\bf h}(p)=K(\cdot,\overline p)$.

By Proposition \ref{prop2.4}, we have
\begin{align}\label{g_rp}
f_{A_n}^{*}(p)=\sum_{j=1}^n\sum_{k=1}^n f(p_j)\widetilde{a_{j,k}^{(n)}}K(p,\overline p_k), \quad p\in E,
\end{align}
where all $p_k$ are distinct points. We need to study what will happen if $p_j=p_k,$ for $j\neq k$.  $(\ref{g_rp})$ depends on the matrix $A_n=(a_{j,k})_{n\times n}$, where $a_{j,k}=K(p_j,\overline p_k)$. Evidently, $A_n$ is singular if $p_j=p_k, j\neq k$. Therefore, $f_{A_n}^*$ is meaningless in such case.

 Now we interpret $f_{A_n}^*$ as follows. Set $k_{p_j}=K(\cdot,\overline p_j)$. Suppose that $\{k_{p_1},...,k_{p_n}\}$ are linearly independent in $H_K(E)$ when $\{p_k\}_{k=1}^{n}$ are distinct from each other. Let $\{\beta_1,...,\beta_n\}$ be the orthogonalization of $\{k_{p_1},...,k_{p_n}\}$ through the Gram-Schmidt (G-S) orthogonalization process. By the G-S orthogonalization process, we know that $\beta_k$ is a linear combination of $\{k_{p_1},...,k_{p_n}\}$. When this relation is combined with the reproducing property (\ref{interpolation}), we have
\begin{equation}\label{f-oth}
f_{A_n}^*=\sum_{k=1}^n\langle f_{A_n}^*, \beta_k \rangle \beta_k=\sum_{k=1}^n\langle f, \beta_k \rangle \beta_k.
\end{equation}

Based on (\ref{f-oth}), we can revise $f_{A_n}^*$ for the situation that $p_j=p_k, j\neq k.$ More precisely, we use $\{\widetilde \beta_1,...,\widetilde\beta_n\}$, the generalization of $\{\beta_1,...,\beta_n\}$ in the above situation, to revise $f_{A_n}^*$. $\{\widetilde \beta_1,...,\widetilde\beta_n\}$ in the contexts of one complex variable and quaternionic variable has been discussed in \cite{QWa1,QSW}. Recently, $\{\widetilde \beta_1,...,\widetilde\beta_n\}$ is formulated in general Hilbert spaces in \cite{Q}. As result, the concept, Complete Dictionary, is deduced in \cite{Q}. The treatment here follows the same line.

For the purpose of convenience, we only interpret $\{\widetilde \beta_1,...,\widetilde\beta_n\}$ for $n=2$. Suppose that $p_1$ is fixed.  We further assume that $k_p(q)=K(q,\overline p)$ is holomorphic in $q$ and anti-holomorphic in p (e.g. the Szeg\"o kernel ). Set $p=p_1+z$, where $z=re^{i\theta}$. By the G-S orthogonalization process, we have
\begin{align*}
\alpha_{\{p_1\}}=\alpha_1=k_{p_1},\quad \alpha_{\{p_1,p_2,...,p_k\}}=\alpha_k= k_{p_k}-\sum_{j=1}^{k-1}\langle k_{p_j},\frac{\alpha_j}{\|\alpha\|}\rangle \frac{\alpha_j}{\|\alpha\|}, k=1,...,n,
\end{align*}
and
\begin{align*}
\{\beta_1,...,\beta_n\}=\{\frac{\alpha_1}{\|\alpha_1\|},...,\frac{\alpha_n}{\|\alpha_n\|}\}.
\end{align*}
Now we consider the following limit
\begin{align}\label{limpp1}
\begin{split}
\lim_{p\to p_1} \beta_{\{p_1,p\}}&=\lim_{p\to p_1}\frac{\alpha_{\{p_1,p\}}}{||\alpha_{\{p_1,p\}}||}\\
&= \lim_{p\to p_1}\frac{\alpha_{\{p_1,p\}}-\alpha_{\{p_1,p_1\}}}{\sqrt{\langle \alpha_{\{p_1,p\}}-\alpha_{\{p_1,p_1\}},\alpha_{\{p_1,p\}}-\alpha_{\{p_1,p_1\}}\rangle}}\\
&= \lim_{p\to p_1}e^{-i\theta}\frac{\frac{\alpha_{\{p_1,p\}}-\alpha_{\{p_1,p_1\}}}{\overline z}}{\sqrt{\langle \frac{\alpha_{\{p_1,p\}}-\alpha_{\{p_1,p_1\}}}{\overline z},\frac{\alpha_{\{p_1,p\}}-\alpha_{\{p_1,p_1\}}}{\overline z}\rangle}}\\
&= e^{-i\theta}\frac{\frac{\mathrm{d} }{\mathrm{d} {\overline p}}\alpha_{\{p_1,p \}}|_{p=p_1}}{||\frac{\mathrm{d} }{\mathrm{d} {\overline p}}\alpha_{\{p_1,p \}}|_{p=p_1}||}\\
&= e^{-i\theta}\frac{\frac{\mathrm{d}}{\mathrm{d} {\overline p}} k_p |_{p=p_1}-\langle \frac{\mathrm{d}}{\mathrm{d} {\overline p}} k_p |_{p=p_1},\frac{\alpha_{1}}{||\alpha_{1}||} \rangle \frac{\alpha_{1}}{||\alpha_{1}||} }{||\frac{\mathrm{d}}{\mathrm{d}{\overline p}} k_p |_{p=p_1}-\langle \frac{\mathrm{d}}{\mathrm{d} {\overline p}} k_p |_{p=p_1},\frac{\alpha_{1}}{||\alpha_{1}||} \rangle \frac{\alpha_{1}}{||\alpha_{1}||}||},
\end{split}
\end{align}
where $\theta$ is the phase of the difference $p-p_1=z$ that keeps to be a constant in the process $p\to p_1.$

 We define $\widetilde \beta_{\{p_1,p_1\}}=\lim_{p\to p_1}\beta_{\{p_1,p\}}$. (\ref{limpp1}) also means that $\widetilde \beta_{\{p_1,p_1\}}$ is the product of $e^{-i\theta}$ and the term generated by involving the derivative $\frac{\mathrm{d} }{\mathrm{d} {\overline p}}k_p |_{p=p_1}$ in the G-S orthogonalization process.

In such situation, we can choose a special direction with $\theta=0$. We can inductively define $\widetilde\beta_k$ by involving $\frac{\mathrm{d}^{m_k}}{\mathrm{d} (\overline p)^{m_k}}k_p|_{p=p_k},$ if necessary, in the G-S orthogonalization process, where $m_k$ is the cardinality of the set $\{j; p_j=p_k, j< k\}$.
 Notice that reproducing kernels in $H^2(\mathbb D)$ and $W(\frac{\pi}{h})$ are anti-holomorphic in the second variable. Now we define a revised $f_{A_n}^*$ as follows.\\
 Let $\{p_k\}_{k=1}^n$ be a sequence of points in $E$, $m_k$ be the cardinality of the set $\{j: p_j=p_k, j< k\}$ for each $p_k$ and
$$
\widetilde K(\cdot,\overline p_k)=\frac{\mathrm d^{m_k}}{\mathrm d (\overline p)^{m_k}}k_p|_{p=p_k},
$$
and $A_n=(a^{(n)}_{j,k})_{n\times n}$, where
\begin{align}\label{new-ajk}
a^{(n)}_{j,k}=\langle \widetilde K(\cdot, \overline p_k), \widetilde K(\cdot, \overline p_j)  \rangle.
\end{align}

Define
 \begin{align}\label{new-aveiro}
\widetilde f_{A_n}^*=\sum_{k=1}^n \langle f, \widetilde \beta_k \rangle\widetilde \beta_k=\sum_{j=1}^n\sum_{k=1}^n \langle f, \widetilde K(\cdot, \overline p_j) \rangle\widetilde{a_{j,k}^{(n)}}\widetilde K(\cdot,\overline p_k),
\end{align}
where $\widetilde{a_{j,k}^{(n)}}$ represents a typical element of $\overline{A^{-1}_n}$. (\ref{new-aveiro}) is well defined because of the following facts:
(1) $\widetilde \beta_k$ is a linear combination of $\{\widetilde K(\cdot,\overline p_1),...,\widetilde K(\cdot, \overline p_k)\}$;
(2)
\begin{align}\label{new-rp}
\langle \widetilde f^*_{A_n}, \widetilde K(\cdot, \overline p_k) \rangle=\langle f, \widetilde K(\cdot, \overline p_k) \rangle=f^{(m_k)}(p_k).
\end{align}
Notice that when all the elements in $\{p_k\}_{k=1}^n$ are distinct with each other, then $\widetilde f_{A_n}^*=f_{A_n}^*$. In such sense  $\widetilde f^*_{A_n}$ generalizes $f^*_{A_n}$. Hereafter, we adopt the same notation $f^*_{A_n}$ for both the distinct and non-distinct cases.

Next we consider
\begin{align}\label{min-new}
p^{*}_{n+1} = \arg \min_{p_{n+1}\in E} ||f-f^{*}_{A_{n+1}}||.
\end{align}
The existence of $p_{n+1}^*$ is evident under the conditions (\ref{new-cond1}) and (\ref{new-cond2}), as given in
\begin{lem}\label{lem_new}
For $f\in H_K(E)$, where $K(\cdot,\overline p)$ is anti-holomorphic in $p$, if
\begin{align}\label{new-cond1}
\lim_{p_{n+1} \to \partial E} \frac{|\langle f, \widetilde K(\cdot,\overline p_{n+1}) \rangle|}{\sqrt{a_{n+1,n+1}}} = 0,
\end{align}
and
\begin{align}\label{new-cond2}
\lim_{p_{n+1} \to \partial E} \frac{|a_{j,n+1}|}{\sqrt{a_{n+1,n+1}}} = 0, j=1,2...n,
\end{align}
then
\[
\lim_{p_{n+1} \to \partial E}||f- f^{*}_{A_{n+1}}|| = ||f- f^{*}_{A_{n}}||,
\]
where $\{p_1,...,p_n\}$ are fixed.
\end{lem}
\begin{proof}
Firstly, let $C_n=(c_1^{(n)}, c_2^{(n)},..., c_n^{(n)})$ and $F_n=(\langle f, \widetilde K(\cdot,\overline p_{1})\rangle, \langle f, \widetilde K(\cdot,\overline p_{2})\rangle,..., \langle f, \widetilde K(\cdot,\overline p_{n})\rangle)$.\\  By $(\ref{new-aveiro})$, we have
\[
{ f}^{*}_{A_n}=\sum_{k=1}^nc_k^{(n)} \widetilde K(\cdot,\overline p_k),
\]
where $C_n=F_n\overline{A_n^{-1}}$.
It is easy to verify that
\begin{align*}
||{f}-{f}^{*}_{A_{n}}||^2
& = ||{f}||^2- \overline{C_n} A_n C_n^{T}.
\end{align*}
Since $A_n$ is a Hermitian and positive definite matrix, we have
$\overline{C_n} A_n C_n^{T}=\overline{F_n} A_n^{-1} A_n \overline{A_n^{-1}}^{T}F_n^{T}=\overline{F_n} A_n^{-1}F_n^{T}$.
So,
\begin{align}
\begin{split}
||{ f}-{ f}^{*}_{A_{n}}||^2 &=||{ f}||^2- \overline{F_n} A_n^{-1}F_n^{T}
\end{split}
\end{align}
and
\begin{align}\label{second1}
\begin{split}
|| f-{ f}^{*}_{A_{n+1}}||^2 &=||{f}||^2- \overline{F_{n+1}} A_{n+1}^{-1}F_{n+1}^{T}.
\end{split}
\end{align}
From now on we denote
\begin{equation*}
A_{n+1} =
\left(
  \begin{array}{ccccc}
    a_{1,1} & \cdots & \cdots & a_{1,n} & a_{1,n+1} \\
    a_{2,1} & \cdots & \cdots  & a_{2,n} & a_{2, n+1} \\
    \vdots & \vdots & \ddots  & \vdots & \vdots \\
    a_{n,1} & \cdots & \cdots & a_{n,n} & a_{n,n+1} \\
    a_{n+1,1} & \cdots & \cdots & a_{n+1,n} & a_{n+1,n+1}
  \end{array}
\right).
\end{equation*}
Secondly, let $adj(A_{n+1})=(A_{j,k})^{T}_{(n+1)\times (n+1)}$ be the adjugate matrix of $A_{n+1}$, where $A_{j,k}$, $(\cdot)^{T}$ and $|\cdot|$ denote the $(j,k)$ cofactor of $A_{n+1}$, the transpose of a matrix and the determinant of a matrix, respectively. Let $(b_{j,k}^{(n)})_{n\times n}$ denote the inverse matrix of $A_n$. Hence,
\begin{align}\label{third3}
(b_{j,k}^{(n+1)})_{(n+1)\times (n+1)} = A^{-1}_{n+1}=\frac{adj(A_{n+1})}{|A_{n+1}|}.
\end{align}
The $(j,k)$ cofactor of $A_{n+1}$ is the product of $(-1)^{j+k}$ and the $(j,k)$ minor of $A_{N+1}$. 
The $(n+1)$-$th$ row expansion of $|A_{n+1}|$ is
\[
|A_{n+1}|=\sum_{j=1}^N a_{n+1,j}A_{n+1,j}+ a_{n+1,n+1} A_{n+1,n+1}.
\]
Let $adj(A_n)=(B_{j,k})^{T}_{n\times n}$, where $B_{j,k}$ denotes the $(j,k)$ cofactor of $A_{n}$.
Let $(a^{(j,k)}_{l,m})_{n \times n}$ denotes the matrix that results from deleting the $j$-$th$ row and the $k$-$th$ column of $A_{n+1}$ and $A^{(j,k)}_{l,m}$ be the $(l,m)$ cofactor of $(a^{(j,k)}_{l,m})_{n \times n}$. For $j,k=1,2,...,n$, it is obvious that
\begin{align}
\begin{split}
a^{(j,k)}_{n,m} &=a_{n+1,m},  m \leq k-1,\\
a^{(j,k)}_{n,m} &=a_{n+1,m+1}, k \leq m \leq n,\\
a^{(j,k)}_{l,n} &=a_{l,n+1},  l \leq j-1,\\
a^{(j,k)}_{l,n} &=a_{l+1,n+1}, j \leq l \leq n ,\\
A^{(j,k)}_{n,n} &= \frac{B_{j,k}}{(-1)^{j+k}}.
\end{split}
\end{align}
For $j=1,2,...,n$,  since $A_{n+1,j}$ is the product of $(-1)^{n+1+j}$ and the $(n+1,j)$ minor of $A_{n+1},$ we have
\[
A_{n+1,j}=(-1)^{n+1+j}\sum_{k=1}^n A^{(n+1,j)}_{k,n} a_{k,n+1}.
\]
Therefore,
\begin{align}\label{third4}
\begin{split}
|A_{n+1}|&=\sum_{j=1}^n\sum_{k=1}^n (-1)^{n+1+j}A^{(n+1,j)}_{k,n}a_{n+1,j}a_{k,n+1}+ a_{n+1,n+1} A_{n+1,n+1}\\
& = a_{n+1,n+1}(|A_n|+\sum_{j=1}^n\sum_{k=1}^n (-1)^{n+1+j}A^{(n+1,j)}_{k,n}\frac{a_{n+1,j}a_{k,n+1}}{a_{n+1,n+1}}).
\end{split}
\end{align}

\noindent For $j,k=1,2,...,n$, we consider the $n$-$th$ row expansion of $A_{j,k}$,
\begin{align*}
A_{j,k} &= (-1)^{j+k}\sum_{m=1}^{n-1} A^{(j,k)}_{n,m} a^{(j,k)}_{n,m} + (-1)^{j+k}A^{(j,k)}_{n,n}a^{(j,k)}_{n,n}\\
& = (-1)^{j+k}\sum_{m=1}^{n-1} A^{(j,k)}_{n,m} a^{(j,k)}_{n,m} + B_{j,k}a_{n+1,n+1}.
\end{align*}
For $m=1,2,...,n-1$, the $(n-1)$-$th$ column expansion of $A^{(j,k)}_{n,m}$ is
\[
A^{(j,k)}_{n,m}=\sum_{l=1}^{n-1} d^{(j,k)}_{m,l}a^{(j,k)}_{l,n},
\]
where $d^{(j,k)}_{m,l}$ depends on the cofactor of the matrix that results from deleting the $n$-$th$ row and the $m$-$th$ column of $(a^{(j,k)}_{l,m})_{n \times n}$.\\
Thus,
\begin{align}\label{third5}
\begin{split}
A_{j,k} &= (-1)^{j+k}\sum_{m=1}^{n-1} A^{(j,k)}_{n,m} a^{(j,k)}_{n,m} + (-1)^{j+k}A^{(j,k)}_{n,n}a^{(j,k)}_{n,n}\\
& = a_{n+1,n+1}(B_{j,k} + \sum_{m=1}^{n-1}\sum_{l=1}^{n-1}{(-1)}^{j+k}d^{(j,k)}_{m,l} \frac{a^{(j,k)}_{l,n} a^{(j,k)}_{n,m}}{a_{n+1,n+1}}).
\end{split}
\end{align}
Therefore, for $j,k=1,2,...,n$, by $(21)$, $(\ref{third3})$, $(\ref{third4})$ and $(\ref{third5})$
\begin{align}\label{third1}
\begin{split}
b_{j,k}^{(n+1)} &= \frac{A_{k,j}}{|A_{n+1}|}\\
& = \frac{(B_{k,j} + \sum_{m=1}^{n-1}\sum_{l=1}^{n-1}{(-1)}^{j+k}d^{(k,j)}_{m,l} \frac{a^{(k,j)}_{l,n} a^{(k,j)}_{n,m}}{a_{n+1,n+1}})}{(|A_n|+\sum_{j=1}^n\sum_{k=1}^n (-1)^{n+1+j} A^{(n+1,j)}_{k,n}\frac{a_{n+1,j}a_{k,n+1}}{a_{n+1,n+1}})}\\
& \to b_{j,k}^{(n)}, \quad as \quad |p_{n+1}| \to \partial E.
\end{split}
\end{align}
Similarly, for $j=1,2,...,n$,
\begin{align*}
A_{n+1,j} &=(-1)^{n+1+j}\sum_{k=1}^n A^{(n+1,j)}_{k,n}a_{k,n+1},\\
b_{j,n+1}^{(n+1)} &=\frac{A_{n+1,j}}{|A_n+1|}\\
&=\frac{(-1)^{n+1+j}\sum_{k=1}^n A^{(n+1,j)}_{k,n}a_{k,n+1}}{a_{n+1,n+1}(|A_n|+\sum_{j=1}^n\sum_{k=1}^n (-1)^{n+1+j} A^{(n+1,j)}_{k,n}\frac{a_{n+1,j}a_{k,n+1}}{a_{n+1,n+1}})}
\end{align*}
and for $k=1,2,...,n$,
\begin{align*}
A_{k,n+1} &= (-1)^{n+1+k}\sum_{j=1}^n A^{(k, n+1)}_{n,j}a_{n+1,j},\\
b_{n+1,k}^{(n+1)} &=\frac{A_{k,n+1}}{|A_n+1|}\\
&=\frac{(-1)^{n+1+k}\sum_{j=1}^n A^{(k, n+1)}_{n,j}a_{n+1,j}}{a_{n+1,n+1}(|A_n|+\sum_{j=1}^n\sum_{k=1}^n (-1)^{n+1+j} A^{(n+1,j)}_{k,n}\frac{a_{n+1,j}a_{k,n+1}}{a_{n+1,n+1}})}
\end{align*}
and
\begin{align*}
A_{n+1,n+1} &=|A_n|,\\
b_{n+1,n+1}^{(n+1)} &=\frac{|A_n|}{|A_{n+1}|}\\
& = \frac{|A_n|}{a_{n+1,n+1}(|A_n|+\sum_{j=1}^n\sum_{k=1}^n (-1)^{n+1+j}A^{(n+1,j)}_{k,n}\frac{a_{n+1,j}a_{k,n+1}}{a_{n+1,n+1}})}.
\end{align*}
Finally, from $(\ref{second1})$,
\begin{align}\label{third2}
\begin{split}
||{ f}-{f}^{*}_{A_{n+1}}||&=||{ f}||^2- \overline{F_{n+1}} A_{n+1}^{-1}F_{n+1}^{T}\\
& = ||{ f}||^2-\sum_{k=1}^{n+1}\langle f,\widetilde K(\cdot,\overline p_k) \rangle\sum_{j=1}^{n+1}\overline{\langle f,\widetilde K(\cdot,\overline p_j) \rangle}b_{j,k}^{(n+1)}\\
& = ||{f}||^2-\sum_{k=1}^{n}\langle f,\widetilde K(\cdot,\overline p_k) \rangle\sum_{j=1}^{n}\overline{\langle f,\widetilde K(\cdot,\overline p_j) \rangle}b_{j,k}^{(n+1)}\\
& - \langle f,\widetilde K(\cdot,\overline p_{n+1}) \rangle\sum_{j=1}^{n+1}\overline{\langle f,\widetilde K(\cdot,\overline p_j) \rangle}b_{j,n+1}^{(n+1)}-\overline{\langle f,\widetilde K(\cdot,\overline p_{n+1}) \rangle}\sum_{k=1}^n \langle f,\widetilde K(\cdot,\overline p_k) \rangle b_{n+1,k}^{(n+1)}.
\end{split}
\end{align}
By $(\ref{third1})$,
when $p_{n+1} \to \partial E$, the second term of the third equality of $(\ref{third2})$ tends to\\
 $\sum_{k=1}^{n}\langle f,\widetilde K(\cdot,\overline p_k) \rangle\sum_{j=1}^{n}\overline{\langle f,\widetilde K(\cdot,\overline p_j) \rangle}b_{j,k}^{(n)}=\overline{F_n} A_n^{-1}F_n^{T} $.
 According to (\ref{new-cond1}) and (\ref{new-cond2}),
 the last two terms of the third equality of $(\ref{third2})$ tend to $0$ when $p_{n+1} \to \partial E$.\\
 Therefore,
\[
\lim_{p_{n+1} \to \partial E}||{ f}-{ f}^{*}_{A_{n+1}}||= ||{f}-{ f}^{*}_{A_{n}}||.
\]
 \end{proof}\\
Given a sequence of points $\{p_k\}_{k=1}^{n+1}$, we call the property
$$
\lim_{p_{n+1} \to \partial E} \frac{|\langle f, \widetilde K(\cdot,\overline p_{n+1}) \rangle|}{||\widetilde K(\cdot, \overline p_{n+1})||} = 0
$$
 the \lq\lq boundary vanishing condition (BVC)."
Since in Lemma \ref{lem_new} we consider the case that $\{p_1,...,p_n\}$ are fixed, the BVC is then reduced to
\begin{align}\label{weak-bvc}
\lim_{p_{n+1} \to \partial E} \frac{|\langle f, K(\cdot,\overline p_{n+1}) \rangle|}{||K(\cdot,\overline p_{n+1})||} = 0.
\end{align}
Note that $p_{n+1}$ must be different from $\{p_1,...,p_n\}$ when $p_{n+1}\to \partial E.$
We call (\ref{weak-bvc}) the weak BVC. Thus, under the assumption of Lemma \ref{lem_new}, the conditions (\ref{new-cond1}) and (\ref{new-cond2}) follows from the weak BVC. In some $H_K(E)$, we can indeed show that the BVC holds.

Under the (weak) BVC assumption we have the selection principle (\ref{min-new}) that implies the convergence of $f_{A_n}^*$, as given in

\begin{thm}\label{AAM-CON}
Suppose that all elements of $\{p_j\}^\infty_{j=1}$ are selected under the principle (\ref{min-new}). For $f\in \mathcal H=H_K(E)$, we have
\begin{align}\label{AAM-convergence}
\lim_{n\to \infty}||f-f^*_{A_n}||=0.
\end{align}
\end{thm}

\begin{proof}
By (\ref{new-aveiro}) and the Riesz-Fischer theorem, there exists $f^*_{A_\infty}\in H_K(E)$ such that
\begin{align}\label{lim_RF}
\lim_{n\to\infty} ||f_{A_n}^*-f^*_{A_\infty}||=0.
\end{align}
Suppose that
 \begin{align}
 g=f-f_{A_\infty}^{*}
 \end{align}
 and
 $$
 ||g||=||f-f_{A_\infty}^*||\neq 0.
 $$
We must have $b\not\in \{p_j\}^\infty_{j=1}$ such that
\begin{align}\label{g_b}
|g(b)|= \delta_0>0.
\end{align}

On one hand,
\begin{align}\label{g}
\delta_0=|g(b)|=| f(b)-f_{A_\infty}^{*}(b)|\leq |f(b)-f_{A_n}^{*}(b)|+|f_{A_\infty}^{*}(b)-f_{A_n}^{*}(b)|.
\end{align}
By $(\ref{lim_RF})$, there exists $N_1$ such that $n>N_1$, the second term of $(\ref{g})$
\[
|f_{A_\infty}^{*}(b)-f_{A_n}^{*}(b)|<\frac{\delta_0}{2}.
\]
Hence,
\[
|f(b)-f_{A_n}^{*}(b)|>\frac{\delta_0}{2}.
\]
On the other hand, combining $b\not \in \{p_j\}^\infty_{j=1}$ and (\ref{new-rp}), we have
\begin{align}
f(b)=f_{A_{n,b}}^{*}(b),
\end{align}
where $A_{n,b}$ is defined by (\ref{new-ajk}) corresponding to $(p_1,...,p_n, b)$.
By $(\ref{min-new})$, we have
\begin{align}
||f||^2-||f_{A_{n,b}}^{*}||^2=||f-f_{A_{n,b}}^{*}||^2\geq ||f-f_{A_{n+1}}^{*}||^2=||f||^2-||f_{A_{n+1}}^{*}||^2.
\end{align}
Hence there exists $N_2$ such that $n>N_2$,
\begin{align}
\begin{split}
|f(b)-f_{A_n}^{*}(b)| &= |f^{*}_{A_{n,b}}(b)-f^{*}_{A_{n}}(b)|\\
& \leq ||f_{A_{n,b}}^{*}-f_{A_{n}}^{*}|| ||K(z,\overline b)||,\\
& \leq L(\sqrt{||f_{A_{n,b}}^{*}||^2-||f_{A_n}^{*}||^2}),\\
& \leq L(\sqrt{||f_{A_{n+1}}^{*}||^2-||f_{A_{n}}^{*}||^2})\\
& \leq L{||f_{A_{n+1}}^{*}-f_{A_{n}}^{*}||}\\
& <\frac{\delta_0}{2},
\end{split}
\end{align}
where $L$ depends on $b$.
If $n>\max\{N_1,N_2\}$, then we arrive a contradiction. This proves the theorem.
\end{proof}

From now on, when performing AMUCD to a particular $H_K(E)$ what we need to verify is the BVC.

In the next section, we are concerned with AMUCD on the RKHSs $H^2(\mathbb D)$ and $W(\frac{\pi}{h})$ that requires verifying their respective BVCs. Indeed, it is not so obvious that the BVC holds in $W(\frac{\pi}{h})$. For $W(\frac{\pi}{h})$, we, instead, verify the weak BVC that turns to be also sufficient.

\begin{remark}

$(1)$ Notice that one can also obtain the result given in Lemma \ref{lem_new} by using formula (\ref{f-oth}). We refer the interested readers to \cite[Section 3]{QSW}.

$(2)$ We also conclude that AMUCD is identical with P-OGA. Here we briefly introduce the idea of P-OGA with the complete dictionary that is the collection of all the directional derivatives of the reproducing kernels of $H_K(E)$. Denote by $\mathcal D=\{k_p,p\in E\}$ the complete dictionary. P-OGA is formulated as follows. For $f\in H_K(E)$,
let \[ f=\sum_{j=1}^{n-1} \langle f, \widetilde\beta_j \rangle \widetilde\beta_j +f_n,\]
where $\{\widetilde \beta_1,..., \widetilde\beta_{n-1}\}$ is the G-S orthogonalization of $\{k_{p_1},...,k_{p_{n-1}}\}$ in the generalized sense and $f_n$ denotes the orthogonal remainder.
We are to choose, for the fixed $k_{p_1},...,k_{p_{n-1}},$ a next dictionary element $k_{p_n}$ to satisfy
\[ |\langle f_n, \widetilde\beta_n\rangle|= \sup \{ |\langle f_n, \widetilde{\beta}_n^\prime\rangle|\ : \ k_p\in {\mathcal D}\},\]
where, with a general testing element $k_p\in {\mathcal D},$  $\{\widetilde \beta_1,...,\widetilde\beta_{n-1}, \widetilde {\beta}_n^\prime\}$ is the G-S orthogonalization of $\{k_{p_1},...,k_{p_{n-1}}, k_p\}.$
In \cite{Q} the convergence of P-OGA is proved, and the convergence rate estimation is obtained. For more details about P-OGA, please see \cite{Q}. In \S 3 we give more details about the relation between AMUCD and P-OGA in $H^2(\mathbb D)$.

$(3)$ For a general RKHS $H_K$ the BVC may not hold. We note that the techniques used in the proofs of Lemma \ref{last-lemma} and Lemma \ref{lem5.3} (see \S 3) do not work for proving the BVC in $W(\frac{\pi}{h})$. In order to implement AMUCD in general $H_K$ spaces other than $H^2(\mathbb D)$, we have to verify the related BVC case by case.

\end{remark}

\section{Applications}

\subsection{Hardy space}
In this section, we are concerned with $H^2(\mathbb D)$.
We say $f\in H^2(\mathbb{D})$, if $f$ is analytic on the open unit disc $\mathbb{D}$ and
 \begin{equation}
\|f\|^{2}_{H^2(\mathbb{D})} =\sup_{0\leq r <1}\frac{1}{2\pi} \int_{0}^{2\pi}|f(re^{it})|^{2}dt < \infty.
\end{equation}
$H^2(\mathbb{D})$ is a RKHS equipped with the inner product
\[
\langle f, g \rangle_{H^2(\mathbb D)} = \frac{1}{2\pi}\int_{0}^{2\pi}f(e^{it})\overline{g(e^{it})} dt, \quad f, g \in H^2(\mathbb D),
\]
where the values of $f(e^{it})$ and $g(e^{it})$ are, respectively, the non-tangential boundary limit functions of $f$ and $g$. Its reproducing kernel is the Szeg\"{o} kernel
\begin{align*}
K_S(z, \overline w)= \frac{1}{1-\overline w z},\quad w,z\in \mathbb D.
\end{align*}
One can immediately obtain the following result from Proposition \ref{prop2.4}.
\begin{thm}\label{thm5.1}
For any $f \in { H^2(\mathbb{D})} $ and distinct points $\{z_j\}_{j=1}^n$ in the unit disc, we have
\begin{align}
f_{A_n}^{*}(z)=  \sum_{j=1}^n \sum_{k =1}^n f(z_j) \widetilde{ a_{j,k}^{(n)}} K_S(z,\overline z_{k}),
\end{align}
where $A_n$ is a matrix with entries $a_{j,k} = K_S(z_{j},\overline z_k )=\frac{1}{1-\overline z_k z_{j}}$,\\
and if $\{z_j\}_{j=1}^{\infty}$ is a uniqueness set of $H^2(\mathbb{D})$,
\begin{align}
f(z)= \lim_{n \to \infty} f_{A_n}^{*}(z)\quad \text{in $H^2(\mathbb D)$}.
\end{align}
\end{thm}
Notice that if $\{z_j\}_{j=1}^{\infty}$ satisfies
$$
\sum_{j=1}^{\infty}(1-|z_j|)=\infty,
$$
then $\{z_j\}_{j=1}^{\infty}$ is a uniqueness set of $H^2(\mathbb{D}).$ The converse result also holds. Those are consequences of the result that zeros of any $f\not\equiv 0 \in H^2(\mathbb D)$ satisfy
$$
\sum_{j=1}^\infty (1-|z_j|) <\infty.
$$

\medskip

Let $\{z_k\}_{k=1}^n$ be a sequence of points in $\mathbb D$, $m_k$ be the cardinality of the set $\{j: z_j=z_k, j< k\}$ for each $z_k$ and
$$
\widetilde K_S(z,\overline z_k)=\frac{\mathrm d^{m_k}}{\mathrm d (\overline z)^{m_k}}K_S(\cdot,\overline z)|_{z=z_k}=\frac{(m_k+1)!z^{m_k}}{(1-\overline z_k z)^{m_k+1}}.
$$ The minimization problem (\ref{min-new}) for $H^2(\mathbb D)$ is stated as follows
\begin{equation}\label{min_pro}
z^{*}_{n+1}:= arg \min_{z_{n+1}\in \mathbb{D}} ||f-f^{*}_{A_{n+1}}||_{H^2(\mathbb{D})}.
\end{equation}

\noindent  As shown by Lemma \ref{lem_new} , to justify the existence of $z^{*}_{n+1}$, we only need to verify the weak BVC in $H^2(\mathbb D)$.
Indeed, we can show the BVC in $H^2(\mathbb D)$.
\begin{lem}\label{last-lemma}
For any $f\in H^2(\mathbb D)$ and any fixed integer $k\geq0$,
\begin{align}\label{last-lem}
\lim_{|z|\to 1^-}\frac{|\langle f, \frac{w^{k}}{(1-\overline z w)^{k+1}} \rangle_{H^2(\mathbb D)} |}{||\frac{w^{k}}{(1-\overline z w)^{k+1}}||_{H^2(\mathbb D)}}=0.
\end{align}
\end{lem}
This result was proved in \cite{WQ}. For the self-containing purpose, we include a proof.\\
\begin{proof}
In fact, (\ref{last-lem}) is a consequence of the following facts.
Since polynomials is dense in $H^2(\mathbb D)$, for any $f\in H^2(\mathbb D)$ and any $\epsilon>0$, there exists $N$ such that
$$
||f-P_N||_{H^2(\mathbb D)}<\frac{\epsilon}{2},
$$
where $P_N$ is a polynomial.
We also know that
$$
||\frac{w^{k}}{(1-\overline z w)^{k+1}}||^2_{H^2(\mathbb D)}=\frac{1}{2\pi}\int_{0}^{2\pi}\left|\frac{e^{ikt}}{(1-\overline z e^{it})^{k+1}}\right|^{2}dt.
$$
Due to Theorem 1.7 in \cite{Z}, we have two constants $c^\prime>0$ and $C^\prime>0$ such that
$$
\frac{c^\prime}{(1-|z|^2)^{2k+1}}\leq ||\frac{w^{k}}{(1-\overline z w)^{k+1}}||^2_{H^2(\mathbb D)}\leq \frac{C^\prime}{(1-|z|^2)^{2k+1}}.
$$
Therefore, when $|z|\to 1^-$,
\begin{align*}
\frac{|\langle f, \frac{w^{k}}{(1-\overline z w)^{k+1}} \rangle_{H^2(\mathbb D)} |}{||\frac{w^{k}}{(1-\overline z w)^{k+1}}||_{H^2(\mathbb D)}} &\leq \frac{|\langle f-P_N, \frac{w^{k}}{(1-\overline z w)^{k+1}} \rangle_{H^2(\mathbb D)} |}{||\frac{w^{k}}{(1-\overline z w)^{k+1}}||_{H^2(\mathbb D)}} +\frac{|\langle P_N, \frac{w^{k}}{(1-\overline z w)^{k+1}} \rangle_{H^2(\mathbb D)} |}{||\frac{w^{k}}{(1-\overline z w)^{k+1}}||_{H^2(\mathbb D)}}\\
& \leq {||f-P_N||_{H^2(\mathbb D)}} + \frac{P_N^{(k)}(z)}{(k+1)!} \frac{(1-|z|^2)^{k+\frac{1}{2}}}{\sqrt{c^\prime}}\\
&\leq \frac{\epsilon}{2}+\frac{\epsilon}{2},
\end{align*}
where $P_N^{(k)}(z)$ is the $k$-th derivative of $P_N(z)$.
\end{proof}

\noindent Notice that (\ref{new-cond2}) in $H^2(\mathbb D)$ is a special case of (\ref{last-lem}). Hence, the existence of $z_{n+1}^*$ follows from (\ref{last-lem}). In the following content, we give more details about the equivalence relation of AMUCD and P-OGA. Then, we can also conclude the convergence of $f_{A_n}^*$ by using the results given in \cite{QWa1}.

We consider the modified Blaschke products corresponding to the sequence $\{z_j\}_{j=1}^{n+1}$
\begin{align*}
B_j(z) &= B_{\{z_1,\ldots,z_j\}}(z)
:= \frac{\sqrt{1-|z_j|^2}}{1-\bar z_j z}\prod_{k = 1}^{j - 1}\frac{z-z_k}{1-\bar z_k z}, j = 1,2,\cdots, n+1,
\end{align*}
where $\{B_1,...,B_{n+1}\}$ is generated by the G-S orthogonalization on $\{\widetilde K_S(\cdot,\overline z_1),..., \widetilde K_S(\cdot,\overline z_{n+1})\}$.

By (\ref{new-aveiro}), we have
\begin{align}\label{fa-afd}
f^{*}_{A_{n+1}}=\sum_{j=1}^{n+1} \langle f^{*}_{A_{n+1}}, B_j\rangle_{H^2(\mathbb{D})}B_j
= \sum_{j=1}^{n+1} \langle f, B_j\rangle_{H^2(\mathbb{D})}B_j.
\end{align}
Then, based on the orthonormal property of $\{B_1,...,B_{n+1}\}$, we have
\begin{align*}
||f-f^{*}_{A_{n+1}}||^2_{H^2(\mathbb{D})}=||f||^2_{H^2(\mathbb{D})}-\sum_{j=1}^n| \langle f, B_j \rangle_{H^2(\mathbb{D})} |^2- |\langle f, B_{n+1} \rangle_{H^2(\mathbb{D})} |^2.
\end{align*}
Hence, $(\ref{min_pro})$ is equal to
\begin{equation}\label{min_pro1}
\begin{split}
z^{*}_{n+1} & := \arg \max_{z_{n+1}\in \mathbb{D}}|\langle f, B_{n+1}\rangle_{H^2(\mathbb{D})} |\\
&= \arg \max_{z_{n+1}\in \mathbb{D}} \sqrt{1-|z_{n+1}|^2} |g_{n+1}(z_{n+1})|,
\end{split}
\end{equation}
where $g_{n+1}(z)=(f- \sum_{k=1}^n \langle f, B_k \rangle_{H^2(\mathbb{D})}B_k)\prod_{l=1}^n \frac{1-\bar z_l z}{z-z_l}$ and $g_{n+1}\in H^2(\mathbb{D})$.
The existence of $z^{*}_{n+1}$ in $(\ref{min_pro1})$ and the convergence of $f_{A_n}^{*}$ have been proved by Qian et al in \cite{QWa1,Q1}. Indeed, the study in this paper is originally motivated by the above observations in $H^2(\mathbb D)$.

\subsection{Paley-Wiener space}
In this section, we consider the following integral transform, for $F\in L^2([\frac{-\pi}{h}, \frac{\pi}{h}]),$ $(h>0),$

\begin{equation}\label{13}
f(z) = \frac{1}{2\pi} \int_{ - \pi/h}^{\pi/h} F(t) e^{ -izt} dt.
\end{equation}
The image space of $(\ref{13})$ is called the Paley-Wiener space $W\left( \frac{\pi}{h} \right)$  comprised of all analytic functions of exponential type  satisfying, for some constant $C$,
$$
|f(z)| \leq C \exp  \left ( \frac{\pi |z|} {h} \right ),\quad |z| \to \infty
$$
and
$$
||f||^2_{L^2(-\infty, \infty)}=\int_{\mathbb R} |f(x)|^2 dx < \infty.
$$
$W\left(\frac{\pi}{h}\right)$ is a RKHS, and its reproducing kernel is
\begin{equation}
K_h(z,\overline{w}) = \frac{\sin \frac{\pi}{h}(z - {\overline w})}{\pi(z - {\overline w})} .
\end{equation}
For more information on properties of $W\left(\frac{\pi}{h}\right)$, please see e.g. \cite{higgins,HS,Y}.\\

Immediately, we have
\begin{thm}\label{thm5.2}
For any $f\in W\left(\frac{\pi}{h}\right)$ and $n$ distinct points $\{z_j\}_{j=1}^n$ in the complex plane, we have
\begin{align}
f_{A_n}^{*}(z)=  \sum_{j=1}^n \sum_{k =1}^n f(z_j) \widetilde{ a_{j,k}^{(n)}} K_h(z,\overline z_{k} ),
\end{align}
where $A_n$ is a matrix with entries $a_{j,k} = K_h(z_{j},\overline z_k )=\frac{\sin \frac{\pi}{h}(z_j-\overline z_k )}{\pi(z_j-\overline z_k )}$,\\
and if $\{z_j\}_{j=1}^{\infty}$ is a uniqueness set of $W\left( \frac{\pi}{h}\right)$, then
\begin{equation}
f(z)= \lim_{n \to \infty} f_{A_n}^{*}(z).
\end{equation}
\end{thm}
Notice that Theorem \ref{thm5.2} is the Shannon sampling theorem if $\{z_j\}_{j=1}^{\infty}$ is replaced by $\{jh\}_{j=-\infty}^{\infty}$. In fact, $\{jh\}_{j=-\infty}^{\infty}$ is a uniqueness set of $W\left( \frac{\pi}{h} \right)$ (see e.g. \cite{Y}).

In fact, we do not know whether the BVC holds in $W(\frac{\pi}{h}).$  However, as mentioned previously, the weak BVC is sufficient for our study.
Without loss of generality, we take $h=1$.
The minimization problem is stated as follows
\begin{align}\label{min_pro2}
\begin{split}
z_{n+1}^{*}:
&= \arg \min_{z_{n+1}\in \mathbb{C}} ||f-f^{*}_{A_{n+1}}||_{L^2(-\infty, \infty)}.
\end{split}
\end{align}

The following lemma shows that the weak BVC holds in $W(\pi)$.
\begin{lem}\label{lem5.3}
If $f\in W(\pi)$,
then
\begin{align}\label{lim2}
\lim_{|z|\to \infty}\frac{|f(z)|}{\sqrt{\frac{\sin \pi(z-\overline z)}{\pi (z-\overline z)}}}= \lim_{|z|\to \infty}\frac{|f(z)|}{\sqrt{\frac{e^{2\pi y}-e^{-2\pi y}}{4\pi y}}}=0.
\end{align}
\end{lem}
\begin{proof}
Since $f\in W(\pi)$, there exists $F\in L^2([-\pi, \pi])$ such that
$$
f(z)=\frac{1}{2\pi}\int_{-\pi}^{\pi}F(t)e^{-i zt}dt, \quad z=x+iy,
$$
Assume $|y|\to \infty$ as $|z|\to \infty$. In this situation, we first prove the above result under the condition $F\in L^p([-\pi, \pi])$, $2<p\leq \infty$. \\
For $q$ satisfying $\frac{1}{p}+\frac{1}{q}=1$, by H\"older's inequality, we have
\begin{align*}
\begin{split}
|f(z)| &\leq \frac{1}{2\pi}\int_{-\pi}^{\pi}|F(t)|e^{y t}dt\\
&\leq \frac{1}{2\pi} (\int_{-\pi}^{\pi}|F(t)|^pdt)^{\frac{1}{p}} (\int_{-\pi}^{\pi}e^{q yt}dt)^{\frac{1}{q}}\\
& \leq \frac{1}{2\pi}(\int_{-\pi}^{\pi}|F(t)|^pdt)^{\frac{1}{p}}(\frac{e^{\pi qy}-e^{-\pi qy}}{qy})^{\frac{1}{q}}.
\end{split}
\end{align*}
Hence,
\begin{align*}
\frac{|f(z)|}{\sqrt{\frac{e^{2\pi y}-e^{-2\pi y}}{4\pi y}}} &\leq \frac{1}{2\pi}(\int_{-\pi}^{\pi}|F(t)|^pdt)^{\frac{1}{p}}\frac{(\frac{e^{\pi qy}-e^{-\pi qy}}{qy})^{\frac{1}{q}}}{\sqrt{\frac{e^{2\pi y}-e^{-2\pi y}}{4\pi y}}}.
\end{align*}
Then, we only need to prove
\begin{align}\label{lim1}
\lim_{|y|\to \infty}\frac{(\frac{e^{\pi qy}-e^{-\pi qy}}{qy})^{\frac{2}{q}}}{{\frac{e^{2\pi y}-e^{-2\pi y}}{4\pi y}}}=0.
\end{align}
When $y\to +\infty$, we have
\begin{align*}
\frac{(\frac{e^{\pi qy}-e^{-\pi qy}}{qy})^{\frac{2}{q}}}{{\frac{e^{2\pi y}-e^{-2\pi y}}{4\pi y}}} &\leq \frac{4\pi y}{(qy)^{\frac{2}{q}}}\frac{(e^{\pi q y})^{\frac{2}{q}}}{e^{2\pi y}-e^{-2\pi y}}\\
& = \frac{4\pi}{(q)^{\frac{2}{q}}y^{\frac{2}{q}-1}}\frac{e^{2\pi y}}{e^{2\pi y}-e^{-2\pi y}}\\
& \to 0,
\end{align*}
where $1\leq q<2$.
Therefore, $(\ref{lim1})$ follows. Similarly, for $y\to -\infty$, we have $(\ref{lim1}).$  This indicates that we have proved $(\ref{lim2})$ under the condition $F\in L^p([-\pi, \pi]), 2<p\leq \infty$.\\
For $F\in L^2([-\pi, \pi])$, we first note that $L^p([-\pi, \pi])(2<p\leq \infty) $ is dense in $L^2([-\pi, \pi])$. In other words, for $F\in L^2([-\pi, \pi])$ and any $\epsilon>0$, there exists $G\in L^p([-\pi, \pi]), 2<p\leq \infty$, such that
$$
||F-G||_{L^2([-\pi, \pi])}<{\pi \epsilon}.
$$
Then, for $F\in L^2([-\pi, \pi])$, when $|y|\to \infty$,
\begin{align*}
\frac{|f(z)|}{\sqrt{\frac{e^{2\pi y}-e^{-2\pi y}}{4\pi y}}} &\leq \frac{\frac{1}{2\pi}\int_{-\pi}^{\pi}|F(t)-G(t)|e^{yt}dt}{\sqrt{\frac{e^{2\pi y}-e^{-2\pi y}}{4\pi y}}} + \frac{\frac{1}{2\pi}\int_{-\pi}^{\pi}|G(t)|e^{yt}dt}{\sqrt{\frac{e^{2\pi y}-e^{-2\pi y}}{4\pi y}}}\\
&\leq \frac{\frac{1}{2\pi}||F-G||_{L^2([-\pi,\pi])}||e^{yt}||_{L^2([-\pi,\pi])}}{||e^{yt}||_{L^2([-\pi,\pi])}} + \frac{\frac{1}{2\pi}||G||_{L^p([-\pi,\pi])}||e^{yt}||_{L^q([-\pi,\pi])}}{||e^{yt}||_{L^2([-\pi,\pi])}}\\
& < \frac{\epsilon}{2}+ \frac{\epsilon}{2},
\end{align*}
where $\frac{1}{p}+\frac{1}{q}=1$.
\medskip
In the next step, we prove that $(\ref{lim2})$ is still true in the situation that $|y|\leq y_0<\infty$ as $|z|\to \infty$. \\
In fact, the dominator of $(\ref{lim2})$ is bounded in this situation. We only need to prove
$$
\lim_{|x|\to \infty}|f(x+iy)|=0.
$$
Since $span\{\frac{\sin{\pi(z-j)}}{\pi(z-j)}; j\in \mathbb Z\}$ is
dense in $W(\pi)$, we simplify the remaining discussion. Specifically, we need to prove, for any fixed $j\in \mathbb
Z$,
\begin{align}\label{lim_p}
\lim_{|x|\to \infty} \left|\frac{\sin{\pi(z-j)}}{\pi(z-j)}\right|
=0.
\end{align}
By directly calculating, $\left|\sin{\pi(z-j)}\right|$ is bounded, when $|y|\leq y_0<\infty$. Then we have (\ref{lim_p}).\\
In the last step, we prove that (\ref{lim2}) without any restriction.
Indeed, we can conclude this by using the previous two conclusions i.e.
\begin{align}\label{case1}
\lim_{|y|\to \infty}\frac{|f(z)|}{\sqrt{\frac{\sin \pi(z-\overline z)}{\pi (z-\overline z)}}}=0,\quad \text{uniformly in } x\in \mathbb R
\end{align}
and
\begin{align}\label{case2}
\lim_{|x|\to\infty} \frac{|f(z)|}{\sqrt{\frac{\sin \pi(z-\overline z)}{\pi (z-\overline z)}}}=0,\quad \text{for $|y|\leq y_0<\infty$}.
\end{align}

We want to show
$$
\lim_{|z|\to\infty} \frac{|f(z)|}{\sqrt{\frac{\sin \pi(z-\overline z)}{\pi (z-\overline z)}}}=0.
$$
Thus means for any given $\epsilon>0$, we need to find $M$ such that for $|z|>M$
$$
\frac{|f(z)|}{\sqrt{\frac{\sin \pi(z-\overline z)}{\pi (z-\overline z)}}}<\epsilon.
$$
By (\ref{case1}) (note that (\ref{case1}) is uniformly in $x$) we can choose $M_1$ such that
for $|y|>M_1$
$$
\frac{|f(z)|}{\sqrt{\frac{\sin \pi(z-\overline z)}{\pi (z-\overline z)}}}<\epsilon.
$$
For $M_1$ being fixed and $|y|\leq M_1$, by (\ref{case2}) we can find $M_2$ large enough such that for $|x|>M_2$
$$
\frac{|f(z)|}{\sqrt{\frac{\sin \pi(z-\overline z)}{\pi (z-\overline z)}}}<\epsilon.
$$
Combining the above facts, we have the following conclusion.
For any $\epsilon>0$, we can find $M^2=M_1^2+M_2^2$  such that for $|z|>M$
$$
\frac{|f(z)|}{\sqrt{\frac{\sin \pi(z-\overline z)}{\pi (z-\overline z)}}}<\epsilon.
$$
Hence
$$
\lim_{|z|\to\infty}\frac{|f(z)|}{\sqrt{\frac{\sin \pi(z-\overline z)}{\pi (z-\overline z)}}}=0.
$$

\end{proof}
\bigskip

\bigskip
\end{document}